\documentclass[11pt,a4paper]{article}

\usepackage{amssymb,amsbsy,amsmath,amsfonts,amssymb,amscd, mathrsfs}
\usepackage{enumerate}

\usepackage{color}
\usepackage[english]{babel}
\usepackage[T1]{fontenc}
\usepackage{indentfirst}
\makeatletter
\@addtoreset{equation}{section}
\makeatother

\newtheorem{statement}{}[section]
\newtheorem{definition}[statement]{Definition}
\newtheorem{theorem}[statement]{Theorem}
\newtheorem{proposition}[statement]{Proposition}

\newtheorem{lemma}[statement]{Lemma}
\newtheorem{rem}{Remark}

\newcommand{\C}{\mathbb C}

\newcommand{\N}{\mathbb N}

\newcommand\e{{\rm e}}
\newcommand{\eps}{\varepsilon}
\newcommand\ind{{\rm 1\kern-.30em I}}

\newcommand\dis{\displaystyle}

\newcommand{\biindice}[3]%
{%

\begin{array}[t]{c}
{\displaystyle #1}\\
{\scriptstyle #2}\\
{\scriptstyle #3}
\end{array}

}

\begin{document}

\title{\bf Essential norms of Volterra and Ces\`aro operators on M\"untz spaces}\date{}
\author{Ihab Al Alam\footnote{Lebanese University, Faculty of Sciences II, Department of Mathematics, P.O. Box 90656 Fanar-Matn, Lebanon,
E-mail: \texttt{ihabalam@yahoo.fr}}, \, Georges Habib\footnote{Lebanese University, Faculty of Sciences II, Department of Mathematics, P.O. Box 90656 Fanar-Matn, Lebanon,
E-mail: \texttt{ghabib@ul.edu.lb}},\, Pascal Lef\`evre \footnote{Laboratoire de Math\'ematiques de Lens (LML), EA 2462, F\'ed\'eration CNRS Nord-Pas-de-Calais FR~2956, Universit\'e d'Artois, rue Jean Souvraz S.P. 18, 62307 Lens Cedex, France, E-mail: \texttt{pascal.lefevre@univ-artois.fr}},
\, Fares Maalouf \footnote{Universit\'e Saint Joseph, Campus des Sciences et Technologies (ESIB), Mar Roukos-Mkall\`es, P.O. Box 1514 Riad El Solh Beirut 1107 2050, Lebanon, E-mail: \texttt{fares.maalouf@usj.edu.lb}}}

\maketitle 

\begin{abstract}
\noindent We study the properties of the Volterra and Ces\`aro operators viewed on the $L^1$-M\"untz space $M_\Lambda^1$ with range in the space of continuous functions. These operators are neither compact nor weakly compact. We estimate how far from being (weakly) compact they are by computing their (generalized) essential norm. It turns out that this latter does not depend on $\Lambda$ and is equal to $1/2$.
\end{abstract}



\section{Introduction}

\noindent Throughout this paper, we denote by $C=C([0,1])$ the space of continuous functions on $[0,1]$ equipped with the supremum norm, and $L^p=L^p([0,1])$ for $p\ge1$, the usual spaces of integrable functions over $[0,1]$ with respect to the Lebesgue measure. \\

\noindent Let $\Lambda$ be a strictly increasing sequence of positive numbers. We say that $\Lambda$ satisfies the {\it M\"untz condition} if and only if $\sum_{\lambda\in\Lambda\setminus \{0\}}1/\lambda<\infty.$ Let $M_{\Lambda}$ be the linear space spanned by the monomials $x^{\lambda}$ with $\lambda\in\Lambda,$ and $M_{\Lambda}^{p}$ the closure of $M_{\Lambda}$ in $L^p$ ($1\leq p<\infty$). The classical theorem of M\"untz states that $M_{\Lambda}^{p}$ is a proper subspace of  $L^{p}$ if and only if the M\"untz condition holds (see \cite{Muntz}, \cite{BorErd}).  This result remains true for the closure of $M_{\Lambda}$ in $C,$ denoted by $M_{\Lambda}^\infty.$ The spaces $M_{\Lambda}^{p}$ ($1\leq p\leq \infty$) are called {\it M\"untz spaces} whenever the M\"untz condition is fulfilled.\\

\noindent An interesting question in this framework is the study of operators acting on M\"untz spaces. In this context, the composition and multiplication operators have been investigated in \cite{A1}, \cite{AL} and \cite{Noor}, where several results have been obtained for the compactness and the weak compactness. However, these operators do not map the M\"untz spaces into themselves in general.\\

\noindent In this paper, we are interested in two classical operators of Hardy-Volterra type, namely the Volterra and Ces\`aro operators. Recall that the Volterra operator $V$ is defined over the space $L^1$ by $V(f)(x)=\displaystyle\int_0^x f(t)\,dt$, and is mapped to the space $C$ by the usual Lebesgue theorem. The restriction $V_\Lambda$ of $V$ to the space $M_\Lambda^1$ does not stabilize it, but almost, since  $V\big(M_\Lambda\big)=M_{1+\Lambda}$. On the other hand, the Ces\`aro operator $\Gamma$ is better-behaved relatively to this problem. Indeed, it is defined  for $x\in (0,1]$ by
$\displaystyle\Gamma(f)(x)=\frac{1}{x}V(f)(x)$ for any $f\in L^1$, and therefore $\Gamma\big(M_\Lambda\big)=M_{\Lambda}.$ Thanks to the Hardy inequality, it is well-known that the Ces\`aro operator maps the space $L^p$ into $L^p$ when $p>1$ but it does not map $L^1$ into itself. By the Clarkson-Erd\"{o}s theorem,  any $f\in M_\Lambda^1$ is equal to an analytic function $\tilde f$ on $[0,1)$, and the a priori ambiguous value $f(0)$ will be defined by taking $f(0):=\tilde f(0)$. Therefore, the Ces\`aro operator on $M_\Lambda^1$ is defined by
$$\Gamma_\Lambda(f)(x)=
\left\{
\begin{array}{ll}
f(0) & \hbox{if $x=0$} \\\\
\dfrac{1}{x}\displaystyle\int_0^xf(t)dt & \hbox{ if $x\in]0,1]$.}
\end{array}
\right.$$

\noindent This can also be written as $\Gamma_\Lambda(f)(x)=\displaystyle\int_0^1f(xu)du$. We easily see that $\Gamma_\Lambda$ maps $M_\Lambda^1$ to $C.$\\

\noindent This paper is divided into four parts. In section 2, we recall some preliminaries on the M\"untz spaces and several notions of operator theory. In section 3, we prove a general criterion  for getting a lower bound for the essential norm of a bounded operator $T$ (see Proposition \ref{PropEss}). Recall that the essential norm of an operator $T$ is the distance from $T$ to the space of compact (or weakly compact) operators. Section 4 will be devoted to the study of some properties of the Volterra and Ces\`aro operators. We show that their essential norm is equal to $\displaystyle \frac{1}{2}$, and is thus independent from the sequence $\Lambda$ (see Theorems \ref{lem:norm} and \ref{thm:ces}).  In the last section, we study some further types of weak compactness, namely the strict singularity and the finite strict singularity; It turns out that $V_\Lambda$ and $\Gamma_\Lambda$ are finitely strictly singular, and that the growth of their Bernstein numbers is of order $\displaystyle \frac{1}{n}$ (Theorem \ref{thm:bernstein}). Actually, in the case of real valued functions, the value of the $n^{\rm th}$-Bernstein number is precisely $1/(2n-1)$ (hence does not depend on $\Lambda$).

\section{Preliminaries}

\noindent We first recall some basic ingredients of the geometry of M\"{u}ntz spaces. These results can be found in detail in \cite{BorErd} and \cite{GuLu}.\\

\noindent Given a strictly increasing sequence of non-negative real numbers $\Lambda=(\lambda_{k})_{k=0}^{\infty}$  such that $\sum_{k=1}^{\infty}1/\lambda_{k}<\infty,$ A. Clarkson and P. Erd\"{o}s (see \cite[p.81]{GuLu} and \cite{Schw59}) proved that any element in the M\"untz spaces can be represented as a series of power functions $x^\lambda$ with $\lambda\in \Lambda.$  The full statement of the result is the following:

\begin{proposition}\label{prop:1}
Let $\Lambda=(\lambda_{k})_{k=0}^{\infty}$ be a strictly increasing sequence of non-negative real numbers satisfying
the M\"untz condition. Assume that the gap condition $\inf\{\lambda_{k+1}-\lambda_{k}:k\in \mathbb N\}>0$ holds. For every function $f\in L^p,$ we have
$f\in M_{\Lambda}^{p}$ if and only if $f$ admits an Erd\"os decomposition, {\it i.e.} for every $x\in[0,1)$, we have
$$f(x)=\sum_{k=0}^{\infty}a_{k}x^{\lambda_{k}}\,.$$
\end{proposition}

\noindent If the gap condition does not hold, any function $f \in M_{\Lambda}^{p}$ can still be  represented by an analytic function on $\{z\in\mathbb{C}\setminus(-\infty,0]: |z|<1\}$ restricted to $(0,1)$. Note that there is no information on the value of $f$  at the point $1$, but one can estimate the supremum norm of $f$ over a compact set far from $1$, by the $L^1$-norm of $f$ over a compact set approaching  $1$ (see \cite[p.185]{BorErd}). Using this last result and the Arzela-Ascoli theorem, the following property of the M\"untz spaces was established in \cite[Corollary 2.5]{AL} :

\begin{lemma} \label{Ihab} Assume that $(f_{n})_{n}$ is a bounded sequence in $M_\Lambda^1$. There is a subsequence of $(f_{n})_{n}$ that converges uniformly on every compact subset of $[0,1)$.
\end{lemma}

\noindent We recall now several notions of operator theory. We begin with the definition of the essential norm of any bounded operator relatively to a closed subspace of operators.

\begin{definition} Let $X$, $Y$ be two Banach spaces, and ${\cal I}$ a closed subspace of operators of  the space $B(X,Y)$  of bounded operators from $X$ to $Y$. The essential norm of $T\in B(X,Y)$, relatively to ${\cal I}$, is the distance from $T$ to ${\cal I}$ defined by
$$\|T\|_{_{e,{\cal I}}}=\inf\{\|T-S\|;\, S\in{\cal I}\}.$$
This is the canonical norm on the quotient space $B(X,Y)/{\cal I}$.
\end{definition}

\noindent With the same notations as in the above definition, if ${\cal I}$ is the space ${\cal K}(X,Y)$ of compact operators from $X$ to $Y$, then $\|T\|_{_{e,{\cal I}}}$ will be  denoted by $\|T\|_{_{e}}$, and will be referred to as \emph{the essential norm of $T$}. If ${\cal I}$ is the space ${\cal W}(X,Y)$ of weakly compact operators from $X$ to $Y$, then $\|T\|_{_{e,{\cal I}}}$ will be denoted by $\|T\|_{_{e,w}}$. The following fact is a direct consequence of the definitions.

\par\ \par \noindent
{\bf Fact:} Let $T: X\to Y$ be an operator, and let $X_0$ and $Y_1$ be subspaces of $X$ and $Y$ respectively, such that $T(X)\subset Y_1$. Let $T_0$  be the operator obtained from $T$ by restricting the domain to $X_0$, and $T_1$  be the one obtained from $T$ by restricting the codomain to $Y_1$. Then \[\|T_0\|_{_{e}} \leq \|T\|_{_{e}}\leq \|T_1\|_{_{e}}.\]

\noindent We recall some other classical notions used in the sequel (see for instance \cite{DJT})

\begin{definition} A Banach space $X$ has the Schur property if every weakly convergent sequence is actually norm convergent.
\end{definition}
It is well known that $\ell^1$ (and consequently all its subspaces) has the Schur property.

\begin{definition}\label{DP}

An operator from a Banach space $X$ to a Banach space $Y$ is weakly compact if $ {T(B_X)}$ is relatively weakly compact.

An operator from a Banach space $X$ to a Banach space $Y$ is Dunford-Pettis (or completely continuous) if it maps any weakly null sequence in $X$ into a norm null sequence in $Y$.
\end{definition} 
Let us mention that every compact operator is Dunford-Pettis. On the other hand any operator from a space $X$ with the Schur property is Dunford-Pettis.\medskip

\noindent Finally, we recall 

\begin{definition}
An operator from a Banach space $X$ to a Banach space $Y$ is nuclear if there exists two sequences  $(y_n)$ in $Y$ and  $(x_n^\ast)$ in $X^\ast$ such that $\sum\|y_n\|.\|x_n^\ast\|<\infty$ and 
$$T(x)=\sum_{n\ge 0}x_n^\ast(x) y_n$$
for every $x\in X$.
\end{definition}
In other words, a nuclear operator is a summable sum of rank one operators. It is {\it a fortiori} a compact operator.

\section{Lower bound of the essential norm}

\noindent In this section, we will prove that the essential norm of any operator can be bounded from below by some term involving the ``height of a discontinuity'' of a function at some point. The idea is based on an argument used in \cite{Z} for the classical essential norm, and in \cite{L2} for the essential norm relatively to the weakly compact operators. For this, we say that a sequence $\dis\big(\tilde x_n\big)_{n\in\mathbb N}$ is a block-subsequence of $\dis\big(x_n\big)_{n\in\mathbb N}$ if there is a sequence of non empty finite subsets of integers  $\dis\big(I_m\big)_{m\in\mathbb N}$ with $\max I_m<\min I_{m+1}$ and $c_i\in[0,1]$ such that
$$\tilde x_m=\dis\sum_{j\in I_m}c_jx_j\quad\hbox{and}\quad \dis\sum_{j\in I_m}c_j=1.$$

\begin{lemma}\label{LemEss}
Let $T: X\rightarrow Y$ be a bounded operator from a Banach space $X$ to another Banach space $Y$. Let $\dis\big(x_n\big)_{n\in\mathbb N}$  be  a normalized sequence in $X$ and $\alpha>0$.
\begin{enumerate}
\item Assume that for any subsequence $\dis\big(x_{\varphi(n)}\big)_{n\in\mathbb N}$ and any $g\in Y$, we have $\dis\overline{\lim}\big\|T\big(x_{\varphi(n)}\big)-g\big\|\ge\alpha.$ Then $\dis\|T\|_{_{e}}\ge\alpha. $
\item Assume that  for any block-subsequence $\dis\big(\tilde x_n\big)_{n\in\mathbb N}$ and any $g\in Y$, we have $\dis\overline{\lim}\big\|T\big(\tilde x_n\big)-g\big\|\ge\alpha.$ Then $\dis\|T\|_{_{e,w}}\ge\alpha$.
\end{enumerate}
\end{lemma}
\noindent {\bf Proof.} We prove only $(2)$ since $(1)$ is similar (and easier). Fix any weakly compact operator $S:X\rightarrow Y$ and consider a normalized sequence $\dis\big(x_n\big)_{n\in\mathbb N}.$ From the weak compactness of $S$, we know that there exists a subsequence of $\dis\big(S\big(x_n\big)\big)_{n\in\mathbb N}$ that converges weakly to some $g\in Y$. By the Banach-Mazur theorem,  there exists a block-subsequence $\dis\big(\tilde x_n\big)_{n\in\mathbb N}$ of $\dis\big( x_n\big)_{n\in\mathbb N}$ with $\dis\big(S\big(\tilde x_n\big)\big)\rightarrow g$ in norm. Note that $\dis\big(\tilde x_n\big)_{n\in\mathbb N}$ lies in the unit ball of $X$.
Hence,  $\dis\|T-S\|\ge\overline{\lim}\big\|T\big(\tilde x_n\big)-g\big\|-\lim\|S\big(\tilde x_n\big)-g\big\|\ge\alpha$. The lemma follows.
\hfill$\square$\\

\noindent Let $K$ be a metric space. We will say that a function $H:K\rightarrow\C$ has a {\sl discontinuity at the point $t_0$ with height $h>0$} if $\dis\overline{\lim_{\delta\rightarrow0}}\dis\, \hbox{diam}\,\big( H(B(t_0,\delta)\big)\ge h$, where $B\big(t_0,\delta\big)$ denotes the open ball centered at $t_0$ with radius $\delta$. Equivalently, one can find two sequences $\dis\big(s_n\big)_{n\in\mathbb N}$ and $\dis\big(s'_n\big)_{n\in\mathbb N}$ both converging to $t_0$ and satisfying $\overline{\lim}\big|H\big(s_n\big)-H\big(s'_n\big)\big|\ge h$.
\begin{proposition}\label{PropEss}
Let $T: X\rightarrow C(K)$ be a bounded operator from a Banach space $X$ to the space $C(K)$ of continuous functions over a compact $K$. We assume that there exists a normalized sequence  $\dis\big(x_n\big)_{n\in\mathbb N}$ in $X$ such that $T\big(x_n\big)$ converges pointwise to some function $H$ with  discontinuity of height $h>0$ at some point $t_0\in X$. Then
$$\|T\|_{_{e,w}}\ge\dis\frac{h}{2}\,\cdot $$
\end{proposition}
\noindent {\bf Proof.} Fix any $g\in C(K)$ and any block-subsequence $\dis\big(\tilde x_n\big)_{n\in\mathbb N}$ of $\dis\big(x_n\big)_{n\in\mathbb N}.$ With the previous notations, we write $\dis\tilde x_m=\dis\sum_{j\in I_m}c_jx_j$, with $\dis\sum_{j\in I_m}c_j=1$. For every $t\in K,$ we have
$$\dis T\big(\tilde x_m\big)(t)=\dis\sum_{j\in I_m}c_jT(x_j)(t)\longrightarrow H(t)$$
since $\min I_m\rightarrow+\infty$. Thus, we obtain that
$$\overline{\lim}\big\|T\big(\tilde x_n\big)-g\big\|\ge\sup_{t\in K}\overline{\lim}\big|T\big(\tilde x_n\big)(t)-g(t)\big|=\sup_{t\in K} \big|H(t)-g(t)\big|.$$
Since the function $H$ has a discontinuity at the point $t_0$ with height $h>0$,  we can find two sequences $s_n\rightarrow t_0$ and $s_n'\rightarrow t_0$ such that   $$\dis\overline{\lim}\big|H\big(s_n\big)-H\big(s'_n\big)\big|\ge h\,.$$ Therefore, we get that
$$\overline{\lim}\big\|T\big(\tilde x_n\big)-g\big\|\ge\max \big\{\overline{\lim}\big|H(s_n)-g(t_0)\big|;\overline{\lim}\big|H(s'_n)-g(t_0)\big| \big\}\ge\frac{h}{2}$$
otherwise we would have
$$\dis h\le\overline{\lim}\big|H\big(s_n\big)-H\big(s'_n\big)\big|\le \overline{\lim}\big|H\big(s_n\big)-g(t_0)\big|+\overline{\lim}\big|g(t_0)-H\big(s'_n\big)\big|<h$$
by continuity of $g$  at point $t_0$. This finishes the proof of the proposition.
\hfill$\square$\\

\noindent\textbf{Application :} Let $K$ be a (metric) connected compact set, $\theta\colon K\rightarrow K$ be a non-constant continuous map, and $T$ be  the composition operator from $C(K)$ into itself defined by $f\in C(K)\mapsto f\circ\theta$. If we consider the sequence  $$x_n(t)=\dfrac{nd(t,\alpha)-1}{nd(t,\alpha)+1}$$ where $\alpha\in \mathrm{Im}\theta$, then

$$H(t)=\left(\lim_{n\to\infty}Tx_n\right)(t)=\left\{
                   \begin{array}{ll}
                     -1, & \hbox{if $t\in \theta^{-1}(\{\alpha\})$;} \\
                     1, & \hbox{if not.}
                   \end{array}
                 \right.$$
The function $H$ has a discontinuity at any point of the boundary of $\theta^{-1}(\{\alpha\})$ (which is non empty by connectedness) with height $2$. We get on one hand that $\dis\|T\|_{_{e}}\ge\dis\|T\|_{_{e,w}}\ge1$.
On the other hand, since $\dis\|T\|_{_{e}}\leq \dis\|T\|\leq1$, we deduce  $\dis\|T\|_{_{e}}=\dis\|T\|_{_{e,w}}=1$.

\section{The Volterra and Ces\`aro operators}

\noindent In this section, we will study some properties of the Volterra and Ces\`aro operators. These operators are neither compact nor weakly compact when restricted to the M\"untz spaces. We will show that their (generalized) essential norm is equal to $\displaystyle\frac{1}{2},$ independently of the choice of the sequence $\Lambda.$\\

\noindent Recall first that the Volterra operator $V$ is defined by $$V(f)=\displaystyle\int_0^x f(t)dt$$ for any $f\in L^1.$ The operator $V$ is clearly well-defined from $L^1$ to $C$, and is bounded with norm $1.$ It also acts as an isometry on the set of positive functions. If considered as an operator  from $L^p$ to $C$ with $p>1$, or from $L^1$ to $L^p$ with $p$ finite, it is easy to see (via Ascoli's theorem) that $V$  is compact. This implies that it is also compact, when viewed as a map from $L^p$ to itself, for any $p\ge1$. Considered as an operator from $L^1$ to $C$, it is easy to see that $V$ is not compact.

Actually the restriction of $V$ to the M\"untz space $M_{\Lambda}^{1}$ is not weakly compact. Indeed, consider any increasing sequence $(\lambda_n)$ in $\Lambda$ and the associated normalized sequence $\dis x\mapsto(\lambda_n+1)\,x^{\lambda_n}$ (for $n\in\N$) in  $M_{\Lambda}^{1}\subset L^1$. The image of this sequence by $V$ is the sequence $\dis x\mapsto x^{\lambda_n+1}$, which admits no weakly convergent subsequence in $C$. We  first compute the essential norm of the Volterra operator viewed as a map from $L^{1}$ to $C.$

\begin{theorem}\label{lem:norm}
Let $V$ be the Volterra operator viewed as a map from $L^{1}$ to $C.$  Then
$$\displaystyle \|V\|_e= \|V\|_{_{e,w}}= \frac{1}{2}\,\cdot$$
\end{theorem}

\noindent {\bf Proof.} For the upper bound, we find a compact operator $S:L^1\rightarrow C$ such that $\displaystyle \|V-S\|\leq \frac{1}{2}.$ For this, we set $\displaystyle S(f)=\frac{1}{2}\displaystyle\int_0^1f(t)dt,$ for  $f\in L^1.$ The operator $S$ is clearly compact since it has rank $1,$ and
\begin{eqnarray*}
\|V-S\|&=& \sup\left\{\|V(f)-S(f)\|_\infty\,:\|f\|_1=1\right\}\\
&=&\sup\left\{\left|\int_0^xf(t)dt-\frac{1}{2}\int_0^1f(t)dt\right|\,: \|f\|_1=1; x\in [0,1]\right\}\\
&=&\sup\left\{\left|\frac{1}{2}\int_0^xf(t)dt-\frac{1}{2}\int_x^1f(t)dt\right|\,: \|f\|_1=1; x\in [0,1]\right\}\\
&\leq &\frac{1}{2}\sup\left\{\int_0^1|f(t)|dt: \|f\|_1=1\right\}=\frac{1}{2}\cdot
\end{eqnarray*}

\noindent For the lower bound, consider any increasing sequence $(\gamma_n)_n$ which tends to infinity and define a normalized sequence $\dis g_n(x)=(\gamma_n+1)\,x^{\gamma_n}$ in $L^1.$ The sequence $V\big(g_n\big)=x^{\gamma_n+1}$ converges pointwise to the function $H$ which is equal to $0$ on $[0,1)$ and to $1$ at the point $1$. The function $H$ has a jump  of height $1$ at $1$, and the result follows by Proposition \ref{PropEss}.
\hfill$\square$\\

Let us mention that the preceding proof actually shows that all the {\it approximation numbers} $a_n(V)$ ({\it i.e.} the distance from $V$ to operators with rank less than $n-1$) are equal to $\frac{1}{2}$ for every $n\ge2$ since we proved indeed
$$\dis\frac{1}{2}\ge a_2(V)\ge\cdots\ge a_n(V)\ge\cdots\ge\|V\|_e\ge\frac{1}{2}\cdot$$ 
We recall that the $a_1(V)$ is equal to $\|V\|=1$.

\medskip

\noindent In the sequel, we will study the properties of the Ces\`aro operator $\Gamma.$ Recall that for any function $f\in L^1,$  $\Gamma(f)$ is defined at any point $x\in (0,1]$ by $$\dis\Gamma(f)(x)=\frac{1}{x}V(f)(x).$$ As  mentioned earlier, the Ces\`aro operator does not map $L^1$ into itself. We will  then study the restriction $\Gamma_\Lambda$  of $\Gamma$ to the M\"untz spaces. Recall that for every $f\in M_\Lambda^1$
$$\Gamma_\Lambda(f)(x)=
\left\{
\begin{array}{ll}
f(0) & \hbox{if $x=0$} \\\\
\dfrac{1}{x}\displaystyle\int_0^xf(t)dt & \hbox{ if $x\in(0,1]$}
\end{array}
\right..
$$

\noindent It is a clear that $\Gamma_\Lambda(f)\in M_\Lambda^\infty$. Moreover, it has the Erd\"os decomposition $\dis\sum_{k=0}^\infty a_k\frac{x^{\lambda_k}}{\lambda_{k}+1}$ when $f$ has the decomposition $\dis\sum_{k=0}^\infty a_k x^{\lambda_k}.$

\begin{proposition} The operator $\Gamma_\Lambda$ satisfies the following:
\begin{enumerate}
\item It is a bounded, one-to-one operator, and  its image $\mathrm{Im}\,\Gamma_\Lambda$ is dense but not closed in $M_\Lambda^\infty$.
\item It is a Dunford-Pettis and non weakly compact operator.
\end{enumerate}
\end{proposition}

\noindent {\bf Proof.} The operator $\Gamma_\Lambda$ is the composition of the bounded operator $V_\Lambda$, and the division operator $Q$ defined from $\displaystyle M_{1+\Lambda}^\infty$ to $\displaystyle  M_{\Lambda}^\infty$ by
\[\displaystyle  Q(f)=\frac{1}{x}f(x).\]
This latter operator is bounded by the closed graph theorem and Proposition \ref{prop:1}.  The operator $\Gamma_\Lambda$ is clearly injective. The density of the image is given by the fact that $M_\Lambda\subset {\rm Im}\,\Gamma_\Lambda.$ Assume for a contradiction that $\mathrm{Im}\,\Gamma_\Lambda$ is closed.  Then $\Gamma_\Lambda:M_\Lambda^1\mapsto M_\Lambda^\infty$ is onto, thus an isomorphism by the Banach theorem. Let $S$ be a lacunary subsequence of $\Lambda$. The restriction of $\Gamma_\Lambda$ to $M_S^1$  is then an isomorphism onto $M_S^\infty$. On the other hand, $M_S^1$ is isomorphic to the space $\ell^1$, and $M_S^\infty$ is isomorphic to the space $ c$ of convergent sequences (see Theorem 9.2.2 of \cite{GuLu}). So $\ell^1$ is isomorphic to $c$, and this is a contradiction since, for instance, the Schur property holds in $\ell^1$ and not in $c$.

\noindent Now we show that $\Gamma_\Lambda$ is not weakly compact. Assume for a contradiction that it is not the case, and consider the normalized sequence $(f_n)_n$ in $M_\Lambda^1$ defined by $f_n(x)=(\lambda_n+1)x^{\lambda_n}$. There exists a subsequence $(f_{n_k})_k$ such that $\Gamma_\Lambda( f_{n_k})=x^{\lambda_{n_k}}$  converges weakly ({\it i.e.} pointwise) to  $f\in C$. Hence, $f=0$ on $[0,1)$ and $f(1)=1$ which is a contradiction. Finally, $\Gamma$ is a Dunford-Pettis operator, because $M_\Lambda^1$ has the Schur property (see the remark after Def.\ref{DP}). Indeed $M_\Lambda^1$ can be realized as an isomorphic copy of a subspace of $\ell^1$ (see \cite{W}).
\hfill$\square$

\begin{rem} We mention the following useful facts:
\begin{enumerate}
\item For any $p\geq 1$, the Ces\`aro operator $\Gamma_\Lambda$ viewed from $M_\Lambda^1$ to $M_\Lambda^p$ is a compact operator. Indeed,  it is the composition of the identity map $i_p$ from $ M_\Lambda^\infty$ to $M_\Lambda^p$, which is compact, with the operator $\Gamma_\Lambda.$ The compactness of $i_p$ is guaranteed by Lemma \ref{Ihab}. Similarly, it follows from the Ascoli theorem that the Ces\`aro operator viewed from $M_\Lambda^p$ to $M_\Lambda^\infty$  is compact  for $p>1$.
\item The operator $D:M_\Lambda^\infty \rightarrow M_\Lambda^1; f\to (xf)'_{|[0,1)}$ is not a bounded operator, since this would give that $\Gamma_\Lambda$ is an isomorphism and $\Gamma_\Lambda^{-1}=D.$
\end{enumerate}
\end{rem}

\begin{theorem} \label{thm:ces} Let $\Lambda$ be an increasing sequence of positive numbers satisfying the  M\"untz condition. Let $\Gamma_\Lambda$ and $V_\Lambda$ be the Ces\`aro and Volterra operators, viewed as maps from $M_\Lambda^1$ to $M_\Lambda^\infty$ or $C$. Then
\begin{enumerate}
\item $\dis\big\|\Gamma_\Lambda\big\|_{e}=\big\|\Gamma_\Lambda\big\|_{_{e,w}}=\frac{1}{2}\,\cdot$
\item $\dis\big\|V_\Lambda\big\|_{e}=\big\|V_\Lambda\big\|_{_{e,w}}=\frac{1}{2}\,\cdot$
\end{enumerate}
\end{theorem}

\noindent To get the upper bound in the both preceding statements, we will use the following proposition, which gives a more general result for some weighted versions of the Volterra-Ces\`aro type operators. For this, fix a continuous function $q$ on $[0,1]$. We define an operator $H_q$ on $M_\Lambda^1$ as follows.  For $f\in M_\Lambda^1$,
 
$$H_q(f)(x)=
\left\{
\begin{array}{ll}
q(0)f(0) & \hbox{if $x=0$} \\\\
\displaystyle \frac{q(x)}{x}\int_0^xf(t)\,dt & \hbox{ if $x\in(0,1]$.}
\end{array}
\right.
$$
If $q(x)=x$, we recover the Volterra operator $V_{\Lambda}$, and if $q=\ind$, we recover the Ces\`aro operator $\Gamma_\Lambda$. We also point out that the range of $H_q$ is included in the space $Q_\Lambda=q.M_\Lambda^\infty.$

\begin{proposition} \label{LemUBess}
With the preceding notations, we have

$$\dis d\Big(H_q,{\cal K}\big(M_\Lambda^1,Q_\Lambda\big)\Big)\le\frac{\|q\|_\infty}{2}\quad\hbox{and}\quad\dis d\Big(H_q,{\cal K}\big(M_\Lambda^1,C\big)\Big)\le\frac{|q(1)|}{2}\,\cdot$$
\end{proposition}

{\noindent \bf Proof.} For the first inequality, let $\lambda\in\Lambda$ and $\rho\in (0,1)$. Let  $R$ and $T_\rho$ be the operators from $\dis M_\Lambda^1$ to $Q_\Lambda$ defined for $f\in M_\Lambda^1$ by
$$\dis R(f)(x)=\frac{q(x)x^\lambda}{2}\int_0^1f(t)\;dt$$
and
\[T_\rho(f)(x)=\left\{\begin{array}{ll}
\rho\, q(0) f(0) & \hbox{if $x=0$} \\\\
\displaystyle \frac{q(x)}{x}\int_0^{\rho x}f(t)\,dt & \hbox{ if $x\in(0,1]$}
\end{array}
\right..
\]
The operator $R$ is compact since it has rank $1$. The operator $T_\rho$ is also compact since it is nuclear. Indeed, it can be written as the sum of an absolutely convergent series of rank one operators in the following way: consider the functionals $\e_n$ defined on $M_{\Lambda}^1$ by
$$\e_n\,\,:\,\,\,\sum_{n\geq 0} a_nx^{\lambda_n}\longmapsto a_n.$$
By a result from \cite[p.178]{BorErd}, the norm of the map $\e_n$ is bounded from above by $C_\rho \rho^{-\lambda_n}$, where $C_{\rho}$ is a constant depending only on $\rho$. Hence for every polynomial $f$ in $M_{\Lambda}^1,$ we have
$$T_\rho(f)(x)=\sum_{n\geq 0} \e_n(f) \frac{\rho^{\lambda_n+1}}{\lambda_n+1} x^{\lambda_n}q(x).$$
The operator $\dis f\longmapsto \e_n(f) \frac{\rho^{\lambda_n+1}}{\lambda_n+1} x^{\lambda_n}q(x)\,$ has rank one, and
$$\sum_{n\geq 0}  \frac{\|\e_n\|\rho^{\lambda_n+1}}{\lambda_n+1}\le\sum_{n\geq 0}  \frac{C_\rho}{\lambda_n+1}<\infty.$$
We show now that $\dis\underline{\lim}_{\rho\rightarrow1^-}\|H_q-(R+T_\rho)\|\le\dis\frac{\|q\|_\infty}{2}.$ For this, fix an arbitrary $\eps\in(0,1)$, and let $c\in(0,1)$ be such that  $\dis 2-c^{\lambda+1}\le(1+\eps)c.$ Let $f$ be in the unit ball of $M_{\Lambda}^1$, and  $x\in(0,1)$.
\begin{eqnarray}\label{eq:1}
\Big(H_q(f)-R(f)-T_\rho(f)\Big)(x)&=&\frac{q(x)}{x}\int_{\rho x}^xf(t)\;dt-\frac{x^\lambda q(x)}{2}\int_0^1f(t)\;dt\nonumber\\
&=&\frac{q(x)}{2}\int_0^1\varphi_x(t)f(t)\;dt\nonumber\\
\end{eqnarray}
with $\dis \varphi_x(t)=\dis\frac{2}{x}-x^\lambda$ if $t\in[\rho x, x]$ and $\dis \varphi_x(t)=-x^\lambda$ if $t\notin[\rho x, x]$. Now we give an upper bound to the value of  $|H_q(f)-R(f)-T_\rho(f)|(x)$  according to the position of $x.$ If $x\leq c$, we apply the first equality in \eqref{eq:1} and obtain the estimate
\begin{eqnarray*}
|H_q(f)(x)-R(f)(x)-T_\rho(f)(x)|&\le & (1-\rho)\|q\|_\infty\sup_{t\in[0,c]}|f(t)|+\frac{\|q\|_\infty}{2}\|f\|_1\\
&\le &\Big((1-\rho)N_\eps+\frac{1}{2}\Big)\|q\|_\infty\|f\|_1\\
&\le&\Big(\eps+\frac{1}{2}\Big)\|q\|_\infty
\end{eqnarray*}
where we used \cite[p.185]{BorErd} to get $N_\eps\ge1$ such that
\[ \sup_{t\in[0,c]}|f(t)|\le N_\eps \|f\|_1\]
for any $f\in M_{\Lambda}^1$, and we chose $\dis\rho$ to be equal to $\displaystyle 1-\frac{\eps}{N_\eps}\in(0,1).$ If $x\ge c$, we use the second equality in  \eqref{eq:1} which yields
$$|H_q(f)(x)-R(f)(x)-T_\rho(f)(x)|\le\frac{1}{2}\|\varphi_x\|_\infty\|q\|_\infty\|f\|_1\le\frac{1+\eps}{2}\|q\|_\infty.$$
Therefore, since $\eps$ is arbitrary, the proof of the first inequality is finished. 

The proof for the second inequality works as well replacing the operator $R$ by
$$\dis R_1(f)(x)=\frac{q(1)}{2}\int_0^1f(t)\;dt.$$
Now the computation follows the lines of the first part but a new function $ \varphi_x$ appears which is given by $\dis \varphi_x(t)=\dis\frac{2q(x)}{x}-q(1)$ if $t\in[\rho x, x]$ and $\dis \varphi_x(t)=-q(1)$ if $t\notin[\rho x, x]$. A suitable choice of $c$ relatively to the continuity of $q$ at the point $1$ gives $\dis\|\varphi_x\|_\infty\le|q(1)|+\eps$ for $x\ge c$. On the other hand, we manage the case $x\le c$ like before to get an upper bound $\dis\eps+\frac{|q(1)|}{2}.$
\hfill$\square$\\

\noindent We note that if $0\in\Lambda$, the previous proposition is not needed for estimating the distance of the Volterra operator $V_\Lambda$ to $\dis{\cal K}\big(M_\Lambda^1,M_\Lambda^\infty\big)$. The short argument of Theorem ~\ref{lem:norm} works for this case.

\medskip

{\noindent \bf Proof of Theorem ~\ref{thm:ces}.} The lower bound follows from Proposition \ref{PropEss}, exactly as in Theorem \ref{lem:norm}. Just pick the exponents $\gamma_n$ in $\Lambda$. The upper bound follows immediately from Proposition ~\ref{LemUBess}.
\hfill$\square$

\section{Strict singularity}

\noindent In this section, we will study other weak forms of compactness, namely the strict singularity and the finite strict singularity. We will show that $V_\Lambda$ and $\Gamma_\Lambda$ share these types of compactness. Moreover, we will estimate the $n^{th}$ Bernstein numbers of these operators by showing that their growth are of order $\dis\frac{1}{n}\,\cdot$ For this, let us recall the two definitions:

\begin{definition}\label{def:ss}
An operator $T$ from a Banach space $X$ to a Banach space $Y$ is strictly singular if it never induces an isomorphism on an infinite dimensional (closed) subspace of $X.$ That is for every $\eps>0$ and every infinite dimensional subspace $E$ of $X,$ there exists $v$ in the unit sphere of $E$ such that $\|T(v)\|\leqslant\eps$.
\end{definition}

\noindent This notion is now very classical and widely studied (see \cite[p.75]{LT} for instance).

\begin{definition}\label{def:fss}
An operator $T$ from a Banach space $X$ to a Banach space $Y$ is finitely strictly singular if for every $\eps>0$, there exists $N_\eps\geqslant1$ such that for every subspace $E$ of $X$ with dimension greater than $N_\eps$, there exists $v$ in the unit sphere of $E$ such that $\|T(v)\|\leqslant\eps$.
\end{definition}

\noindent This latter definition can be reformulated in terms of the so-called {\it Bernstein approximation numbers} (see \cite{Pl} for instance). Recall that the $n^{th}$ Bernstein number of an operator $T$ is defined as
$$\dis b_n(T)=\biindice{\sup}{E\subset X}{dim(E)=n}\biindice{\inf}{v\in E}{\|v\|=1}\|T(v)\|.$$
Hence, with this terminology, the operator $T$ is finitely strictly singular if and only if $\dis\big(b_n(T)\big)_{n\geqslant1}$ belongs to the space $c_0$ of null sequences. This notion has appeared in the late sixties. For instance, in a paper of V. Milman \cite{Mi}, it is proved that the identity from $\ell^p$ to $\ell^q$ ($p<q$) is finitely strictly singular (see \cite{CFPTT}, \cite{Pl}, \cite{L}, \cite{LR} for recent results). It is also well-known that

\medskip 
\centerline{compactness\quad$\Longrightarrow$\quad finite strict singularity \quad$\Longrightarrow$\quad  strict singularity}

\noindent and that the reverse implications are not true. Moreover, complete continuity is not comparable to finite strict singularity in general. In \cite{L}, it is proved that the classical Volterra operator is finitely strictly singular (hence strictly singular) from $L^1$ to $C$ and moreover the Bernstein numbers satisfy
$$b_n(V)\approx\frac{1}{n}\,\cdot$$ By restriction, the operator $V_\Lambda$ is clearly finitely strictly singular as well. It is easy to see that the same property occurs for $\Gamma_\Lambda.$ But a natural question then arises: how fast the Bernstein numbers vanish? Surprisingly, it turns out that it does not really depend on $\Lambda$ and the answer is given in the following theorem. In the statement below, we are interested in $V_\Lambda$ and $\Gamma_\Lambda$ acting from  $M_{\Lambda}^1$ to $C$, but actually it is worth pointing out that it does not change the values of the Bernstein numbers if we consider these operators from $M_{\Lambda}^1$ onto their range.

\begin{theorem} \label{thm:bernstein}
The  operators $V_\Lambda$ and $\Gamma_\Lambda$ are finitely strictly singular. Moreover the growth of their Bernstein numbers is of order $\dis\frac{1}{n}\cdot$ For every $n\ge1$, we have
$$b_{n}\big(V_\Lambda\big)\ge\dis\frac{1}{2n-1}\quad\hbox{and}\quad b_{n}\big(\Gamma_\Lambda\big)\ge\dis\frac{1}{2n-1}\,\cdot$$
In particular, in the case of real valued functions, we have $\dis b_{n}\big(V_\Lambda\big)=\dis\frac{1}{2n-1}\,\cdot$
\end{theorem}

{\noindent \bf Proof.} Since $V$ is finitely strictly singular (\cite{L}), it is clear by restriction that  $V_\Lambda$ is also finitely strictly singular. Moreover by the fact that the map $f\in M_{1+\Lambda}^\infty\longmapsto \dis\frac{1}{x}f\in M_{\Lambda}^\infty$ is bounded, the operator $\Gamma_\Lambda$ is again finitely strictly singular. Moreover, we have that $\dis b_{n}\big(V_\Lambda\big)\le b_{n}\big(V\big)\le\dis C_n$ where $C_n$ is equal to $\frac{1}{2n-1}$ or $\frac{\sqrt2}{n}$ according to the fact that we consider real or complex valued functions (\cite{L}). Now, let us estimate the lower bound. We mainly follow the ideas of Newman in \cite[Lemma 2]{Ne}. Fix $n\ge1$ and $\eps>0$ and let $\big(\lambda'_n\big)_{n\in\mathbb N}$ be a subsequence of $\Lambda$ of positive numbers going very fast to the infinity such that it satisfies the following condition
$$\dis\prod\Bigg[1-\frac{2\lambda_j'^2\big(1+\ln(\lambda_{j+1}')\big)}{\lambda_{j+1}'}\Bigg]\ge1-\eps.$$
It is straightforward from \cite[Lemma 2]{Ne} to get for every $a_1,\ldots,a_n\in\C$ that
$$\Big\|\dis\sum_{k=1}^na_kx^{\lambda'_k}\Big\|_\infty\ge(1-\eps)\max_{1\le m\le n}\Big\|\dis\sum_{k=1}^ma_kx^{\lambda'_k}\Big\|_\infty\ge(1-\eps)\max_{1\le m\le n}\Big|\dis\sum_{k=1}^ma_k\Big|.$$
Next, we consider the space $E$ spanned by $\dis x^{\lambda'_1},\ldots,x^{\lambda'_n}$. Hence, we may write any $f\in E$ as $\dis\sum_{k=1}^na_k\big(\lambda'_k+1\big)x^{\lambda'_k}$ and  we have $\dis \Gamma_\Lambda(f)=\dis\sum_{k=1}^na_kx^{\lambda'_k}.$ Denoting $a=(a_k)_{k=1,\ldots,n},$ and $s_m=\dis\sum_{k=1}^ma_k$ (where $1\le m\le n$), we point out that 
$$\dis\|a\|_{\ell^1}=\big|s_1\big|+\sum_{m=2}^{n}\big|s_m-s_{m-1}\big|\le(2n-1)\max_{1\le m\le n}\big|s_m\big|.$$
Therefore, we find
$$\|\Gamma_\Lambda(f)\|_\infty\ge(1-\eps)\max_{1\le m\le n}\big|s_m\big|\ge\frac{(1-\eps)}{2n-1}\|a\|_{\ell^1}.$$
Finally, the inequality $\|a\|_{\ell^1}\ge\|f\|_1$ finishes the proof of the theorem. The same argument works as well for the operator $V_\Lambda$.
\hfill$\square$

\begin{rem} We point out that the lower bound for the Bernstein numbers of the Volterra operator $V_\Lambda$ turns out to be the same one than the lower bound for $V$ (see \cite{L}). Hence the speed of the Bernstein numbers is independent of $\Lambda$ (up to uniform constants for complex valued functions).
\end{rem}

\noindent {\bf Acknowledgment.} This work was made with the support of the PHC C\`edre project EsFo. The third author would like to thank the colleagues for the warm atmosphere during the stays at the Lebanese University in Beirut.





\begin{thebibliography}{99}
\bibitem[A]{A1} I. Al Alam, {\it Essential norms of weighted composition operators on M\"untz spaces}, J. Math. Anal. {\bf 358} (2009),  273-280.
\bibitem[AL]{AL} I. Al Alam, P. Lef\`evre, {\it Essential norms of weighted composition operators on $L^1$- M\"untz spaces}, Serdica Math. J. {\bf 40} (2014), 241-260.
\bibitem[BE]{BorErd} P. Borwein, T. Erd\'{e}lyi, {\it Polynomials and polynomial inequalities}, Springer, Berlin Heidelberg New York, (1995).
\bibitem[CE]{ClarErd} A. Clarkson, P. Erd\"{o}s, {\it Approximation by polynomials}, Duke Math. J. {\bf 10} (1943), 5-11.
\bibitem[CFPTT]{CFPTT} I. Chalendar, E. Fricain, A. Popov, D. Timotin and V. Troitsky, {\it Finitely strictly singular operators between James spaces},
J. Funct. Anal. {\bf 256} (2009), 1258-1268.
\bibitem[DJT]{DJT} J. Diestel, H. Jarchow, A. Tonge, {\it Absolutely summing operators}, Cambridge University Press (1995).
\bibitem[E]{Erdélyi} T. Erdélyi, {\it The ``full Clarkson-Erdös-Schwartz Theorem'' on the closure of non-dense Müntz spaces},  Studia Math. {\bf 155} (2003), 145-152.
\bibitem[GL]{GuLu} V. I. Gurariy, W. Lusky, {\it Geometry of M\"{u}ntz spaces and related questions}, Lecture Notes in Mathematics, Springer, Berlin Heidelberg New York, (2005).
\bibitem[L]{L} P. Lef\`evre, {\it The Volterra operator is finitely strictly singular from $L^1$ to $L^\infty$}, J. Approx. Theory {\bf 214} (2017), 1-8.
\bibitem[L2]{L2} P. Lef\`evre, {\it Generalized essential norm of weighted composition operators on some uniform algebras of analytic functions}, Integral Equations and Operator Theory {\bf 63} (2009), 557-569.
\bibitem[LR]{LR} P. Lef\`evre, L. Rodr\'iguez-Piazza, {\it  Finitely strictly singular operators in harmonic analysis and function theory},
Advances in Math. {\bf 255} (2014), 119-152.
\bibitem[LT]{LT} J. Lindenstrauss, L. Tzafriri, {\it Classical Banach spaces. I. Sequence spaces}, Vol. 92, Springer-Verlag, Berlin-New York, (1977).
\bibitem[Mi]{Mi} V. Milman, {\it Operators of class $C_{0}$ and $C\sp*_{0}$. (Russian) Teor. Funkcii Funkcional, }
Anal. i Priložen. {\bf 10} (1970), 15-26.
\bibitem[M]{Muntz} Ch. H. M\"{u}ntz, {\it \"{U}ber den Approximationssatz von Weierstrass}, H. A. Schwarz's Festschrift, Berlin, pp. 303-312, (1914).
\bibitem[N]{Noor} S. Waleed Noor, {\it Embeddings of M\"untz spaces: Composition operators}, Int. Equa. and Oper. Theo. {\bf 73} (2012), 589-602.
\bibitem[Ne]{Ne} D. Newman, {\it A M\"untz space having no complement}, J. Approx. Theory {\bf 40} (1984), 351-354.
\bibitem[Pl]{Pl} A. Plichko, {\it Superstrictly singular and superstrictly cosingular operators}, Functional analysis and its applications, North-Holland Math. Stud., 197, Elsevier, Amsterdam, (2004).
\bibitem[S]{Schw59} L. Schwartz, {\it Etude des sommes d'exponentielles}, Hermann, paris, (1959).
\bibitem[W]{W} D. Werner, {\it  A remark about Müntz spaces},  http ://page.mi.fu-berlin.de/werner/preprints/muentz.pdf.
\bibitem[Z]{Z} L. Zheng, {\it The essential norms and spectra of composition operators on $H^\infty$}, Pacific J. Math. {\bf 203}, 503-510 (2002).
\end{thebibliography}
\end{document}